\documentclass[pamm,a4paper,fleqn]{w-art}
\usepackage{times,cite,w-thm}
\usepackage[utf8]{inputenc}

\usepackage[english]{babel}
\usepackage{amsthm}
\usepackage{dsfont}
\usepackage{amssymb}
\usepackage{hyperref}
\theoremstyle{plain}
\usepackage{caption}

\usepackage{tikz}

\newcommand{\R}{\mathbb{R}}

\newcommand{\cd}{\dot{c}}
\newcommand{\chid}{\dot{\chi}}
\newcommand{\taud}{\dot{\tau}}
\newcommand{\hd}{\dot{h}}
\newcommand{\h}{\textup{h}}
\usepackage{graphicx}
\usepackage{animate}
\usepackage{graphicx}
\begin{document}

\TitleLanguage[EN]
\title[Cell seeding dynamics for meniscus regeneration]{Cell seeding dynamics in a porous scaffold material designed for meniscus tissue regeneration
}

\author{\firstname{Henry} \lastname{Jäger}\inst{1,}%
\footnote{Corresponding author: henry.jaeger@edu.rptu.de 
    } 
}

\author{\firstname{Elise} \lastname{Grosjean}\inst{2,}%
     \footnote{elise.grosjean@inria.fr}}


\author{\firstname{Steffen} \lastname{Plunder}\inst{3,}%
     \footnote{plunder.steffen.2a@kyoto-u.ac.jp}}

\author{\firstname{Claudia} \lastname{Redenbach}\inst{1,}
		\footnote{claudia.redenbach@rptu.de}}

\author{\firstname{Alex} \lastname{Keilmann}\inst{1,}
		\footnote{keilmann@rptu.de}}

\author{\firstname{Bernd} \lastname{Simeon}\inst{1,}%
     \footnote{bernd.simeon@math.rptu.de}}

\author{\firstname{Christina} \lastname{Surulescu}\inst{1,}%
     \footnote{surulescu@mathematik.uni-kl.de}}
 
\address[\inst{1}]{\CountryCode[DE]
	Rheinland-Pfälzische Technische Universität Kaiserslautern-Landau (RPTU),
	FB Mathematik, 67663 Kaiserslautern, Germany} 

\address[\inst{2}]{\CountryCode[FR]
	INRIA-Saclay m3disim,
	1 Rue Honor\'e d'Estienne d'Orves, 91120 Palaiseau, France}

\address[\inst{3}]{\CountryCode[JP]
	Institute for the Advanced Study of Human Biology (ASHBi), KUIAS, Kyoto University, Faculty of Medicine Bldg. B, Kyoto, 606-8303, Japan}

\AbstractLanguage[EN]
\begin{abstract}
We study the dynamics of a seeding experiment where
a fibrous scaffold material is colonized by two types of cell populations. The specific application that we have in mind is related to the idea of meniscus tissue regeneration.
In order to support
the development of a promising replacement material, we discuss certain 
rate equations for the densities of human mesenchymal stem cells and chondrocytes and for
the production of collagen-containing extracellular matrix. For qualitative studies,
we start with a system of ordinary differential equations and refine then the model to include 
spatial effects of the underlying nonwoven scaffold structure.
Numerical experiments as well as a complete set of parameters for future benchmarking are provided.
\end{abstract}
\maketitle                   

\section{Introduction}
Meniscal lesions are a frequent injury of the knee joint 
and involve a substantial risk for premature osteoarthritis, in particular if a tear occurs in the
inner, avascular zone. One focus in the search for an alternative treatment lies
in regenerative approaches, and we report here on such research and the corresponding
development of an in-silico experimental environment. More specifically, we study the
situation where an artificial scaffold is seeded with mesenchymal stem cells that differentiate into chondrocytes and build up a stable collagen-containing extracelluar matrix. 

A detailed mathematical model of the involved processes has been proposed 
in \cite{bib2}. It uses a multiscale approach for the chemical, topological and mechanical influences and corresponding processes of the seeding experiment. It starts at the microscopic level of single cells and associated receptor binding dynamics, then passes through the mesoscale of cell distribution functions described by kinetic transport equations, and finally
applies parabolic upscaling to obtain effective equations for the dynamics of 
the macroscopic population densities. However, the resulting diffusion-reaction-advection equations with mechanical coupling terms are quite involved and expensive to handle numerically. In the present work, we introduce a model based on ordinary differential equations (ODEs) as a complement to the partial differential equation (PDE) model and investigate 
relevant data and parameters for future benchmarking. 

The paper is organized as follows. We outline the rate equations of the ODE model in Section 2, extend the description to the PDE model in Section 3 and present numerical experiments in Section 4, along with a complete set of parameters and data.
\section{The rate equations}
Consider an experimental set-up for a porous medium that is to be seeded by human mesenchymal stem cells (hMSCs). These are supposed to differentiate into chondrocytes. This process is controlled by a differentiation medium and also enhanced by a coating of the scaffold with hyaluron. The chondrocytes produce  extracellular matrix (ECM) material that, jointly with the spreading of the cells, gradually fills the pores in the scaffold. 

To describe this seeding experiment,
we introduce the five time-dependent variables $c_1$: representing the density of hMSCs, $c_2$: the density of chondrocytes, $\chi$:  the concentration of the differentiation medium, $h$: the concentration of hyaluron, and $\tau$: the density of the produced ECM. As derived in detail in \cite{bib2}, the following system of rate equations expresses the corresponding cell 
dynamics if spatial effects are omitted. I.e., for the moment we consider all quantities to be homogeneously distributed in space and neglect the motility of the cells that would result in additional taxis and diffusion type terms. 
Under these assumptions, the dominating dynamics can be summarized as
		\begin{align}
			\cd_{1}
			&=-\alpha_1({ \chi},S)c_{1}+\alpha_2({ \chi},S)\frac{\omega_1}{\omega_2}c_{2} +\beta c_{1}\left(1-\frac{c_{1}}{C_1^*}-\frac{c_{2}}{C_2^*}   \right ), \\
			\cd_{2}&=\alpha_1({ \chi},S)\frac{\omega_2}{\omega_1}c_{1}-\alpha_2({ \chi},S)c_{2}, \\
			\chid &=  -a_\chi \left( \frac{c_{1}}{C_1^*} + \frac{c_{2}}{C_2^*} \right) \chi, \label{eq:3}\\ 
			\hd & =  -\gamma_1 \frac{c_{1}}{C_1^*} h -\gamma_2 \frac{c_{2}}{C_2^*} h
			+\gamma_3 \frac{c_2}{1 +\frac{c_2}{C_2^*}}, \label{eq:4} \\
			\taud& = -\delta_1 \frac{c_{1}}{C_1^*} \tau +  \delta_2 c_2\,.\label{eq:5}
		\end{align}
The right hand side for the hMSCs $c_1$ consists of several terms. The first one is a conversion term with a rate $\alpha_1$ that depends on the available differentiation medium $\chi$ and
the mechanical stimulus $S$ that we describe below. This term describes the differentiation of hMSCs into chondrocytes. The second one describes the reciprocal dedifferentiation of chondrocytes into hMSCs with a further rate $\alpha_2$, while the third one models the growth of $c_1$ using a growth rate $\beta$. The growth is slowed down if the pores are getting filled with new cells.
Here, $C_1^*$ and $C_2^*$ are constants which are assumed to represent the carrying capacities of hMSCs and chondrocytes, respectively. 

The growth of chondrocytes $c_2$ is governed by similar terms for differentiation and dedifferentiation, but with opposite signs. The ratios $\omega_1 / \omega_2$ and $\omega_2 / \omega_1 $ stem from the different velocities of the two cell types. In \eqref{eq:3}, 
the differentiation medium $\chi$ is uptaken by both cell types with a constant rate $a_{\chi}$. 
Analogously, the hyaluron is uptaken by both cell types, which is described in the fourth equation with the rates $\gamma_1$ and $\gamma_2.$ The last term in \eqref{eq:4} represents a very limited expression of hyaluron by chondrocytes. Finally, the ECM of the density $\tau$ is produced by chondrocytes and uptaken by hMSCs, which is described in \eqref{eq:5} with the rates $\delta_2$ and $\delta_1$.

\paragraph*{The nonlinear rate functions $\alpha_1$ and $\alpha_2$.}
We make the ansatz $\alpha_1(S,\chi) = \alpha_{1,S}(S) \cdot \alpha_{1,\chi}(\chi)$ and use the positive constants $S_{\textup{min}}$, $S_{\textup{max}}$, $\alpha_{1,\textup{min}}$, $\alpha_{1,\textup{max}}$ and $\chi_c$. For simpler notation, we define $S_\textup{d} := \frac{1}{10}(S_{\textup{max}}-S_{\textup{min}})$ and then choose
	$$\alpha_{1,S}(S) =     
	\begin{cases}
		\alpha_{1,\textup{min}}, & S \leq S_{\textup{min}} - S_\textup{d} \\
		\frac{ 	\alpha_{1,\textup{min}} - 	\alpha_{1,\textup{max}}}{4} \cdot \left( \frac{S- S_{\textup{min}} }{S_\textup{d}} \right) ^3 + 3 \cdot \frac{	\alpha_{1,\textup{max}} - 	\alpha_{1,\textup{min}} }{4} \cdot \left( \frac{S- S_{\textup{min}}}{S_\textup{d}} \right) + \frac{	\alpha_{1,\textup{max}} + 	\alpha_{1,\textup{min}} }{2} , &  S_{\textup{min}} - S_\textup{d} < S \leq S_{\textup{min}} + S_\textup{d} \\
		\alpha_{1,\textup{max}}, &  S_{\textup{min}} + S_\textup{d} < S \leq S_{\textup{max}} - S_\textup{d} \\
		\frac{ 	\alpha_{1,\textup{max}} - 	\alpha_{1,\textup{min}}}{4} \cdot \left( \frac{S- S_{\textup{max}} }{S_\textup{d}} \right) ^3 + 3 \cdot \frac{	\alpha_{1,\textup{min}} - 	\alpha_{1,\textup{max}} }{4} \cdot \left( \frac{S- S_{\textup{max}}}{S_\textup{d}} \right) + \frac{	\alpha_{1,\textup{max}} + 	\alpha_{1,\textup{min}} }{2} , &  S_{\textup{max}} - S_\textup{d} < S \leq S_{\textup{max}} + S_\textup{d} \\
		\alpha_{1,\textup{min}}, &   S_{\textup{max}} + S_\textup{d} < S \\
	\end{cases} $$ as well as
	$$\alpha_{1,\chi}(\chi) =  \frac{\chi^2}{\chi_c^2 + \chi ^2}.$$
	The choice of $\alpha_{1,S}$ yields a $C^1$-function that approximates a piecewise constant function by means of Hermite interpolation.
	For $\alpha_2(S,\chi),$ we analogously make the ansatz $\alpha_2(S,\chi) = \alpha_{2,S}(S) \cdot \alpha_{2,\chi}(\chi)$ and define  
	$$\alpha_{2,S}(S) = 
	\begin{cases}
		\alpha_{2,\textup{max}}, & S \leq S_{\textup{min}} \\
		\alpha_{2,\textup{max}} \cdot \frac{S_{\textup{min}}}{S},
		& S_{\textup{min}} < S \\
	\end{cases}
	$$
	as well as
	$$ \alpha_{2,\chi}(\chi) =  \frac{\chi_c^2}{\chi_c^2 + \chi ^2} $$
	with constants $S_{\textup{min}}$, $\chi_c$ and $\alpha_{2,\textup{max}} .$
In Fig.~\ref{fig:alpha}, the different functions are displayed, with the parameters given in Table  \ref{tab:parameters} below.
 
\begin{figure}
		\includegraphics[width=0.9 \textwidth]{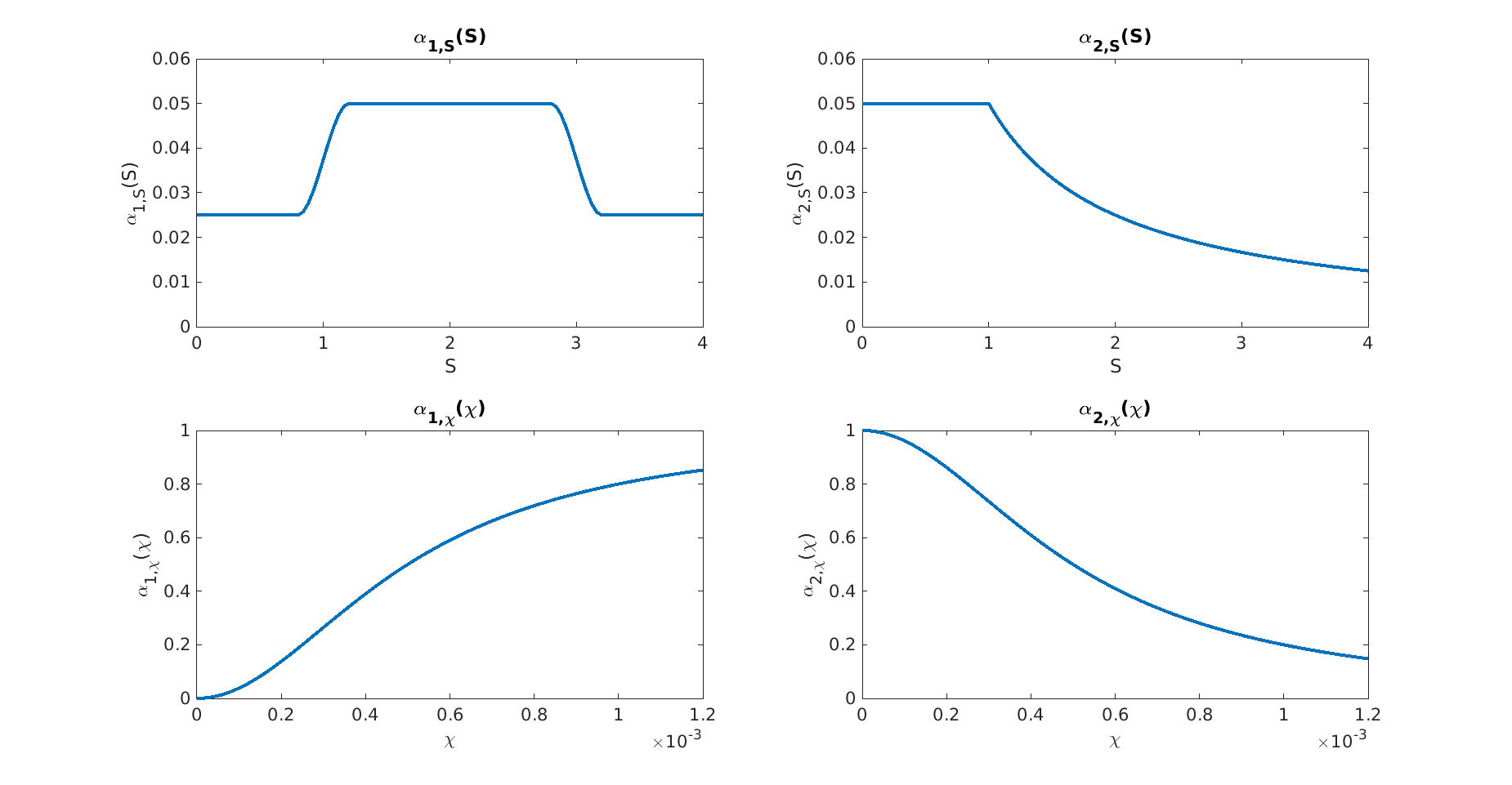}
		\caption{Ansatz functions for $\alpha$ with the parameters from the Table \ref{tab:parameters}}
		\label{fig:alpha}
\end{figure}

%
\section{Background and extension to spatial effects}
In \cite{bib2}, a multiscale modeling procedure for the chemical, topological and mechanical influences and corresponding processes of the seeding experiment is presented. It starts at the microscopic level of single cells and associated receptor binding dynamics, then passes through the mesoscale of cell distribution functions depending on time, position and velocity, and finally
applies parabolic upscaling to obtain effective equations for the dynamics of 
the macroscopic population densities.

The resulting system of reaction-diffusion-advection equations uses the same variables
$c_1, c_2, \chi, h,  \tau$ as the ODE model above, but now all quantities depend on time
$t \in [t_0,t_{\mbox{\tiny end}}]$ and space $x \in \Omega \subset \mathbb{R}^d$ where
$\Omega$ is the domain occupied by the scaffold in $d=2$ or $d=3$ dimensions.
The system reads 
\begin{align}
			 \partial_tc_{1}-\nabla\nabla :\left (\mathbb D_1c_{1}\right )  \nonumber 
			+ & \nabla \cdot \left (\frac{k^-\lambda _{11}}{B(h,\tau)^2(B(h,\tau)+\lambda_{10})}\mathbb D_1\nabla B(h,\tau)c_{1}\right ) \\
			=&-\alpha_1({ \chi},S)c_{1}+\alpha_2({ \chi},S)\frac{\omega_1}{\omega_2}c_{2} +\beta c_{1}\left(1-\frac{c_{1}}{C_1^*}-\frac{c_{2}}{C_2^*}  \right),
			\label{pde-c1} \\
			\partial_tc_{2}-\nabla\nabla :\left (\mathbb D_2c_{2}\right )=&\alpha_1({ \chi},S)\frac{\omega_2}{\omega_1}c_{1}-\alpha_2({ \chi},S)c_{2}, \\
			\partial _t\chi =& D_\chi\Delta \chi -a_\chi \left(
			\frac{c_{1}}{C_1^*} + \frac{c_{2}}{C_2^*} \right)
			\chi, \\ 
			\partial_th = & -\gamma_1 \frac{c_{1}}{C_1^*} h -\gamma_2  \frac{c_{2}}{C_2^*} h
			+\gamma_3  \frac{ c_{2} }{1 + c_{2}/C_2^* },  \\
			\partial_t\tau = & -\delta_1  \frac{c_{1}}{C_1^*}   \tau + \delta_2 c_2\,.
			\label{pde-tau}
		\end{align}
Here, the differential operators on the left hand sides of the equations for $c_1$ and $c_2$
model diffusion and taxis processes by taking into account the fibre orientation distribution in the nonwoven scaffold. More specifically,  the cell diffusion tensors $\mathbb D_1 \in \R^{3 \times 3}$
and $\mathbb D_2 \in \R^{3 \times 3}$ are computed from
the mesoscopic orientation distribution of the fibers using scaffold data
obtained from an a priori imaging analysis. 
The latter is based on evaluating a symmetric and positive-definite parameter matrix $A \in \mathbb{R}^{3  \times 3}$ 
based on an angular central Gaussian (ACG) distribution. For the scaffold that we consider, $A$ can be assumed to be diagonal.
Following the method in \cite{ospald}, this
 leads to the steps
\begin{equation}\label{elliptic-integr-D}
\mathbb D_\beta =c_{A,\beta}\int \limits _0^\infty \prod \limits _{i=1}^3(b_i+\zeta )^{-\frac{\beta _i+1}{2}}d\zeta , \qquad
c_{A,\beta}:=\frac{|\det A|^{-1/2}}{4}\prod \limits_{i=1}^3\frac{\beta_i!}{(\beta_i/2)!},
\end{equation}
where $\beta =(\beta _1,\beta_2,\beta _3)^T$ is a multiindex with $\sum\limits _{i=1}^3\beta _i=2$ used to specify the entries $\mathbb D_\beta$ and $b_i$ ($i=1,2,3$) are the eigenvalues of $A^{-1}$.
The entries $\mathbb D_\beta$ can thus be obtained by computing the elliptic integrals in \eqref{elliptic-integr-D} above.
The $(i, j)$-entry of $\mathbb D_1$ is then given by $\mathbb D_{\beta}$ where $\beta$ is the unique multi-index
such that $\{ k : \beta_k \neq 0 \} = \{ i, j \}$,
multiplied by a prefactor of $ s_1^2 / \lambda_{10}$.
For instance, $\mathbb D_{(2,0,0)}=\frac{|\det A|^{-1/2}}{2}\int \limits _0^\infty(b_1+\zeta )^{-3/2}(b_2+\zeta)^{-1/2}(b_3+\zeta)^{-1/2}d\zeta $ determines the $(1,1)$-entry of $\mathbb D_1$, $\mathbb D_{(1,1,0)}$ the $(1,2)$-entry as well as the $(2,1)$-entry of $\mathbb D_1$,
$\mathbb D_{(1,0,1)}$ the $(1,3)$- and the $(3,1)$- entry of $\mathbb D_1$ and so forth.  Once $\mathbb D_1$
has been computed, the second tensor follows from 
$\mathbb D_2={ \frac{\lambda_{10}}{\lambda_2}\left (\frac{s_2}{s_1}\right )^2\mathbb D_1}$ where 
the positive constants $s_1$ and $s_2$ are the speeds of hMSCs and chondrocytes,
respectively, and $\lambda_{10}, \lambda_2$ are also positive constants, relating to the turning frequency of the two cell phenotypes within the scaffold.

Moreover, in (\ref{pde-c1}) on the left hand side it holds
$B(h,\tau)=k_{1}^{+} \frac{h}{H}+k_{2}^{+} \frac{\tau}{K}+k^-$
where $k_j^+$ and $k_j^-$ denote attachment and respectively detachment rates of hMSC to hyaluron ($j=1$) and ECM ($j=2$) with corresponding reference densities 
$H$ and $K$.
For simplicity we assume $k_1^-=k_2^-=:k^-$. The parameter $\lambda_{11}>0$ is, too, a cell turning rate and stems from the mesoscopic description of the velocity jump process for hMSCs. For more details see \cite{bib2} and references therein.		
Besides initial conditions that are specified in the next section, the PDE system 
(\ref{pde-c1})-(\ref{pde-tau}) requires boundary conditions.
They read
	\begin{align}
		\left (\mathbb D_1\nabla c_1+\left(\nabla \cdot \mathbb D_1
		-\frac{k^-\lambda _{11}}{B(h,\tau)^2(B(h,\tau)+\lambda_{10})}\mathbb D_1\nabla B(h,\tau)\right )c_{1}\right )\cdot \nu &=0\quad \text{on}\quad \partial \Omega, \\
		\left (\nabla \cdot \mathbb D_2c_2+\mathbb D_2\nabla c_2\right )\cdot \nu &=0\quad \text{on}\quad \partial \Omega, \\
		\nabla \chi \cdot \nu  &=0\quad \text{on} \quad \partial \Omega
	\end{align}
and represent no-flux conditions along the boundary for $c_1, c_2$ and $\chi$.

\section {Numerical results}
In this section, we discuss simulation results for both the ODE model and the enhanced PDE model in 2D.
	\subsection{Simulating the ODE model}
	As initial values we set 
\begin{equation}
		c_1(0) = 0.001, \quad
		c_2(0)  = 0, \quad
		 \chi(0) = 0.001, \quad
		h(0) = 1000, \quad 
		\tau(0) = 0. 
	\end{equation}
	The time $t$ is measured in hours, the densities $c_1$ and $c_2$ are measured in $1/\mu m^2$ and the concentration $h$, the concentration~$\chi$ and the density $\tau$ in $\textup{mol}/\mu m^2$.
	For the mechanical stimulus $S,$ we use $S(t) = 0.5 +\textup{cos}(t/10)$.
	The parameters and their units are listed in Table \ref{tab:parameters} below. The numerical time integration has been carried out using MATLAB's ode23s solver \cite{matlab}.

\begin{figure}
		\includegraphics[width=1.1 \textwidth]{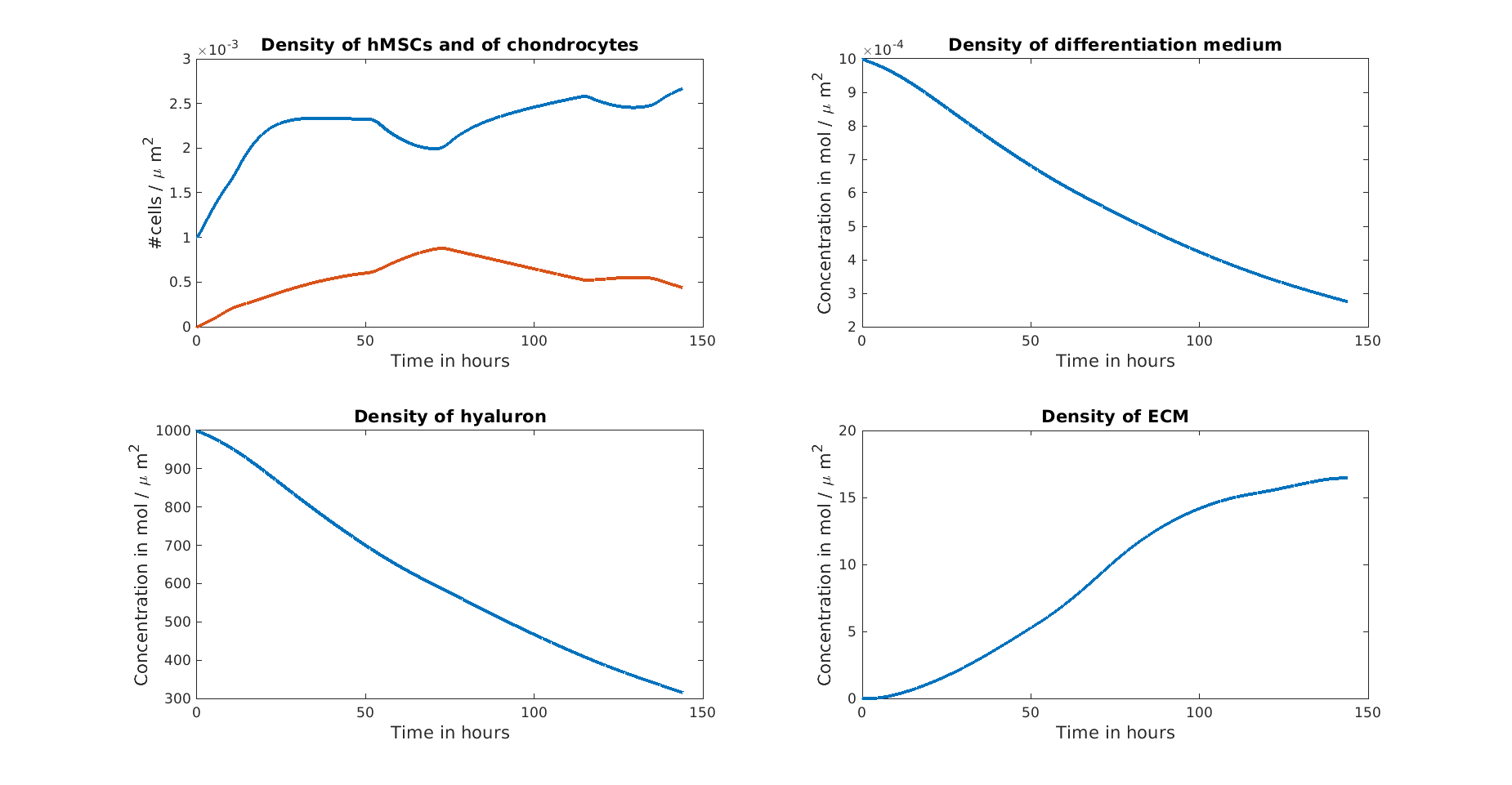}
		\caption{ODE model, $t_{\mbox{{\tiny end}}} =144 \h$. On the top left, the blue curve stands for the hMSCs and the red one for the chondrocytes.}
		\label{fig:ODE1}
	\end{figure}

\begin{figure}
		\includegraphics[width=0.9 \textwidth]{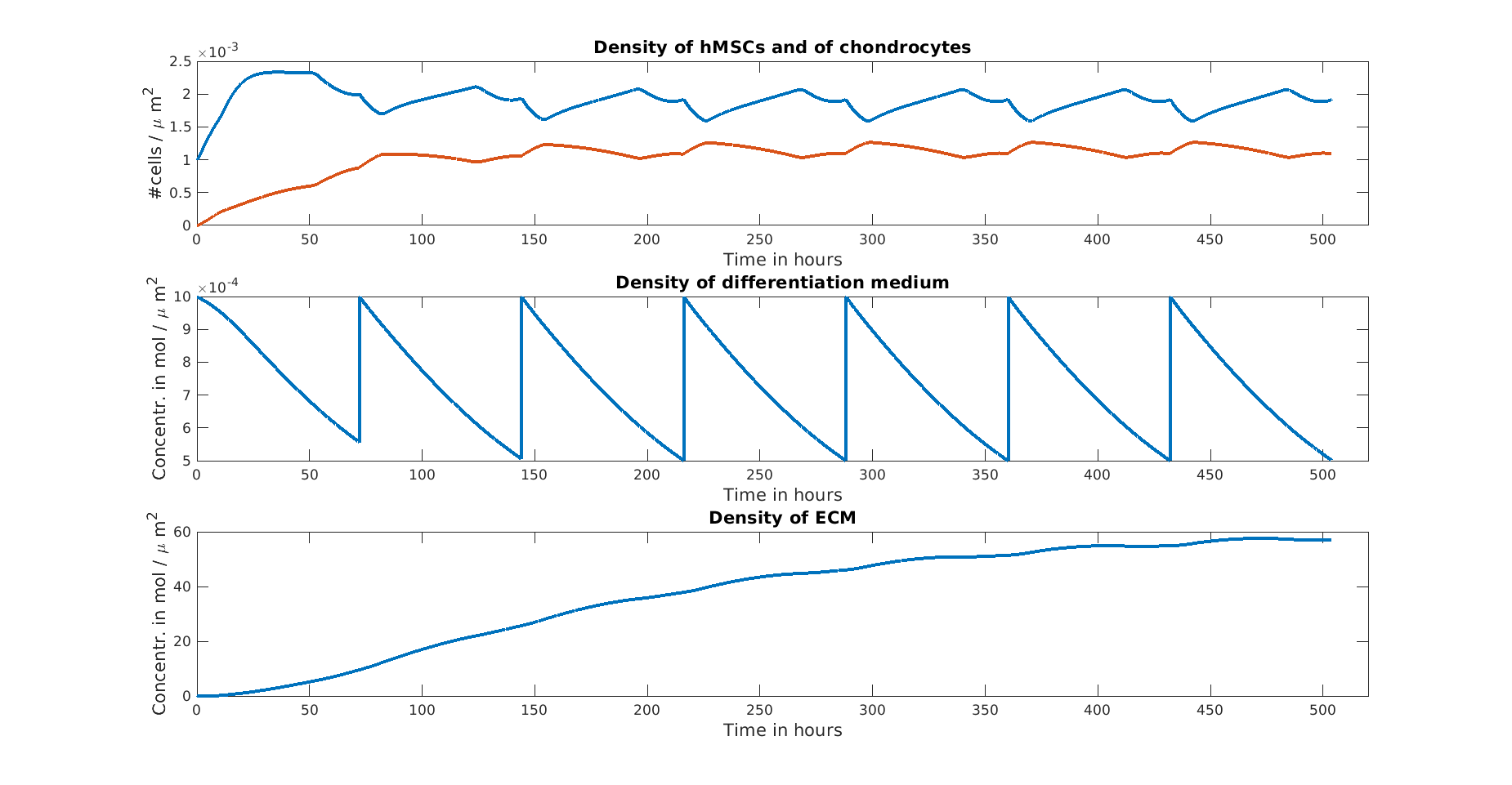}
		\caption{ODE model, evolution of cell densities, differentiation medium and produced ECM over time,
		$t_{\mbox{{\tiny end}}}=504 \h $.}
		\label{fig:ODE2}
	\end{figure}

	 The first in-silico experiment runs over $6$ days or $t_{\mbox{\tiny end}}=144$ hours, respectively,
and the resulting cell dynamics is shown in Figure \ref{fig:ODE1}. 
	One observes that the consumption of the differentiation medium leads to enhanced dedifferentiation of chondrocytes after day $3$, which also slows down the desired  production of collagen-containing ECM.
This effect is also observed in the experiment. As a remedy, the differentiation medium is renewed every $3$ days, and the resulting dynamics shows a clear improvement when performing a longer integration over $21$ days - as Fig.~\ref{fig:ODE2} demonstrates.

	\subsection{Simulating the PDE model}
	The PDE model is solved on a triangular mesh using FreeFEM++ \cite{freefem} with a 
	first order Non-symmetric Interior Penalty discontinuous Galerkin (NIP dG) scheme in space and implicit Euler as robust discretization in time with stepsize $\delta t = 0.1 \h$.
	The circular domain mimics the flat scaffold disks that are used in the 
	real experiment and is given by
	$		\Omega = \{ x \in \R^2: ||x - (2500,2500)^T ||_2 \leq 2500  \}	 $
	(in micrometers).  
	 	As initial values, we used
	\begin{align*}
		c_1(x,0) & = 0.001 \cdot \textup{exp}(-15 \cdot ((x-2500 \mu m)/1000 \mu m)^2 - 15 \cdot ((y-2500 \mu m)/1000 \mu m)^2 ), \\
		c_2(x,0) & = 0,  \quad\chi(x,0) = 0.001, \quad
		h(x,0) = 995 + r, \quad
		\tau(x,0) = 0,
	\end{align*}
	where the five quantities are measured in the same units as in the ODE model and
	where $r$ is an $U(0,1)$-distributed random variable, cf.~the simulations presented
     in \cite{bib1} for a standard diffusion process.
	For simplicity, the term\\
	$k^-\lambda _{11}/({B(h,\tau)^2(B(h,\tau)+\lambda_{10})}) \mathbb D_1$ was replaced by an identity matrix in the simulation.
	
\begin{figure}\centering
		\includegraphics[width =0.45 \textwidth]{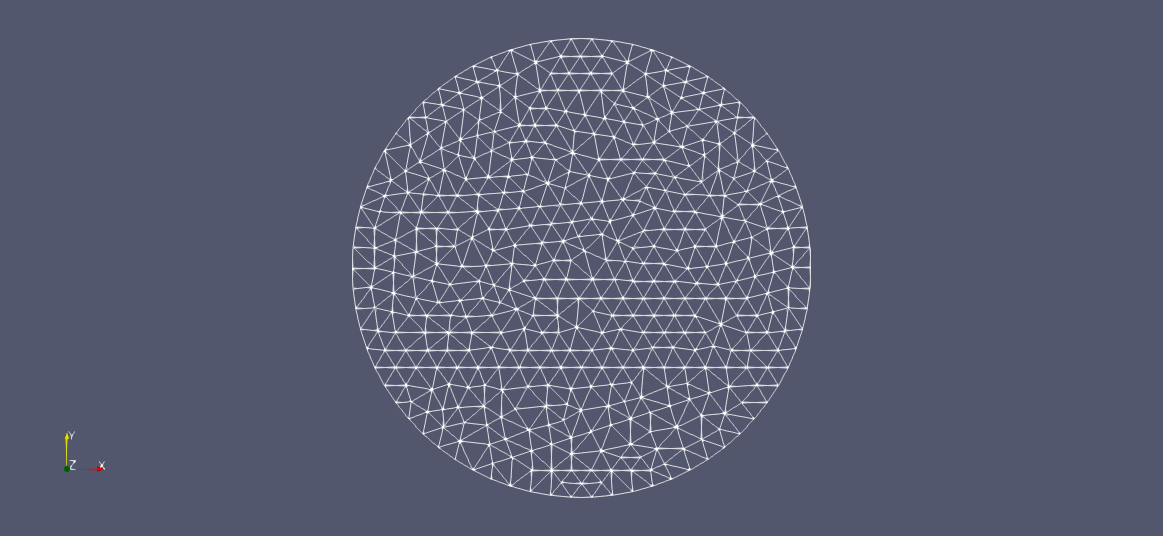}
		\qquad
		\includegraphics[width=0.45 \textwidth]{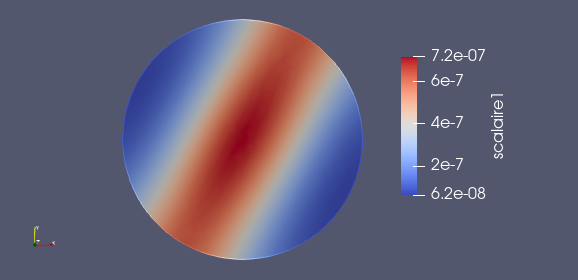}
		 \caption{FEM mesh (left) and snapshot of $c_2$ at $t= 2 \h$ showing the spreading of chondrocytes along the dominating fiber orientation.}
		 \label{fig:PDE1}
	\end{figure}

\begin{figure}
		\includegraphics[width=1.1 \textwidth]{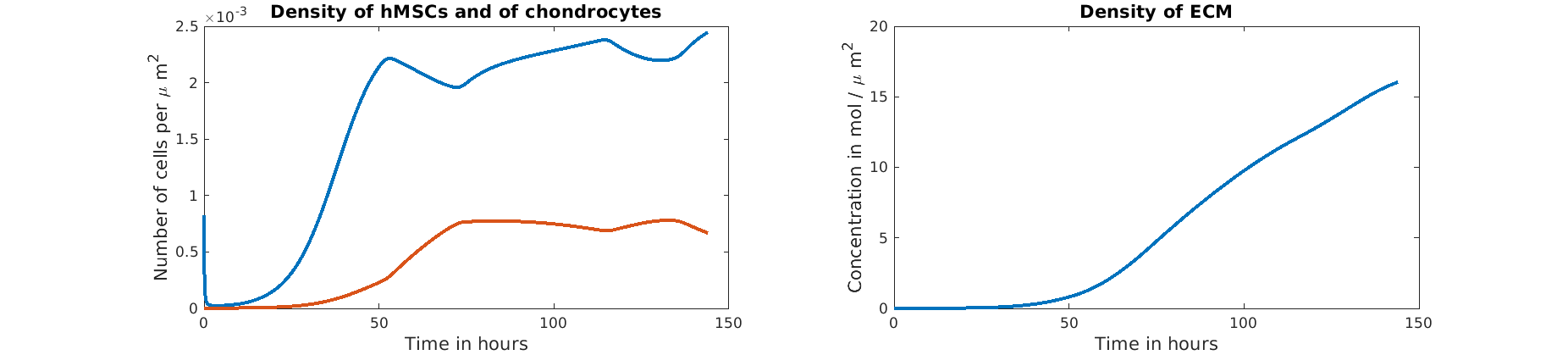}
		\caption{PDE model, evolution of the variables over the time at the midpoint. On the left, the blue curve stands for the hMSCs and the red one for the chondrocytes. Note the peak of $c_1$ at $t=0$ that stems from the initial condition. Due to diffusion, this peak first flattens before the cell dynamics comes into play.}
		\label{fig:PDE2}
	\end{figure}

As Fig.~\ref{fig:PDE1} shows, the diffusion-dominated spread of cells in the scaffold proceeds primarily along a 
diagonal band, which corresponds to the dominating fiber direction in the analyzed samples of nonwoven. Looking at
the temporal behavior of the cell densities in the midpoint of the scaffold, 
we observe in Fig.~\ref{fig:PDE2} a similar pattern as in 
the ODE model (cf. Fig.~\ref{fig:ODE1}), after the initialization phase has passed. Nevertheless, the PDE model predicts a slightly higher level of chondrocyte density at larger times, which is probably due to the tactic effects. However, there is less ECM expression, the onset of which is delayed by the spatial cell spread dominating the early dynamics.
	
\begin{table}\caption{List of the parameter values that have been used in the numerical simulations. Due to the fact that $S$ is only considered in terms of ratios, we handle it as a dimensionless quantity. }\label{tab:parameters}
\begin{tabular}{ccc}
	\hspace{-0.25cm}
\begin{tabular}{|c|c|c|c|}
			\hline
			Parameter & Unit & Value \\[1ex]
			\hline
			$\beta$ & $1/\h$ & 0.5/3 \\[1ex]
			\hline
			$s_1$ & $\mu m/\h $ & 30 \\[1ex]
			\hline
			$s_2$ &  $\mu m/\h $ & 12 \\[1ex]
			\hline
			$\omega_1$ & $(\mu m/\h)^{d-1}$ & 30 \\[1ex]
			\hline
			$\omega_2$ &  $(\mu m/\h)^{d-1}$ & 12 \\[1ex]
			\hline
			$\delta_1$ & $1/\h$ & 3.3 \\[1ex]
			\hline
			$\delta_2$ & $\text{mol}/\h$ & 330 \\[1ex]
			\hline
			$S_{\textup{min}}$ & 1 & 1  \\[1ex]
			\hline
			$S_{\textup{max}}$ & 1 & 3   \\[1ex]
			\hline
			$\alpha_{1,\textup{min}}$ & $1/\h$ & 0.025 \\[1ex]
			\hline
			$\alpha_{1,\textup{max}}$  & $1/\h$ & 0.05 \\[1ex]
			\hline
			$\alpha_{2,\textup{max}}$ & $1/\h$ & $0.05$ \\[1ex]
			\hline
		\end{tabular} \hspace{-0.25cm} 
	&

	\begin{tabular}{|c|c|c|}
		\hline
		Parameter & Unit & Value \\[1ex]
		\hline
		$\alpha_{ \chi}$  & $1/\h$ & 3.18 \\[1ex]
		\hline\texttt{}
		$\gamma_1$ & $1/\h$ & 3.3 \\[1ex]
		\hline
		$\gamma_2$ & $1/\h$ & 1 \\[1ex]
		\hline
		$\gamma_3$ & $1/\h$ & $3.307 \cdot 10^{-3}$  \\[1ex]
		\hline
		$D_{\chi}$ & $\mu m^d/\h$ & $10^6$ \\[1ex]
		\hline
		$k_1^+/H$ & $\mu m^d/(\h\text{ mol})$ & 5 \\[1ex]
		\hline
		$k_2^+/K$ & $\mu m^d/\h$ & 1 \\[1ex]
		\hline
		$C_1^*$ & $ 1 / (\mu m^d)$ & $3.024 \cdot 10^{-3}$ \\[1ex]
		\hline
		$C_2^*$ & $ 1 / (\mu m^d)$ & $3.024 \cdot 10^{-3}$ \\[1ex]
		\hline
		$\lambda_{10}$ & $1/ \h $ & $9 \cdot 10^{-4}$ \\[1ex]
		\hline
		$\lambda_{2}$ & $1 / \h $ & $1.44 \cdot 10^{-4}$ \\[1ex]
		\hline
	\end{tabular} \hspace{-0.5cm}
	&
	\begin{tabular}{l}
		$\mathbb D_1 
		= \frac{s_1^2}{\lambda_{10}} \left( \begin{array}{ccc}
			0.204 & 0.189 & 0.169 \\
			0.189 & 0.447 & 0.251 \\
			0.169 & 0.251 & 0.349 \\
		\end{array}
		\right) $ 
	\end{tabular}
\end{tabular}
\end{table}

\section{Conclusion}

Our in-silico experiment for the cell seeding dynamics in a porous scaffold material 
takes the dominant processes into account. While the ODE model allows a straightforward adaptation of parameters and introduction of further effects, 
the more involved PDE model takes the scaffold 
structure as a nonwoven into account and provides detailed spatial patterns. 
In this way, our work can be expected to form a valuable tool for future in-vitro
experiments designed for meniscus tissue regeneration.

\begin{acknowledgement}
We would like to thank our partners Andreas Seitz and Graciosa Texeira from the
Institute of Orthopaedic Research and Biomechanics at the 
University Ulm as well as 
Martin Dauner, Michael Doser, Carsten Linti and G\"unter Schmidt from the
Deutsche Institute f\"ur Textil- und Faserforschung (DITF) in Denkendorf for 
the fruitful collaboration.
  We are grateful to DFG for funding this work within the Priority Program SPP2311
  \cite{spp2311}.
\end{acknowledgement}

\vspace{\baselineskip}

\end{document}